\documentclass[12pt]{amsart}
\usepackage{amssymb}
\usepackage{amsthm}
\usepackage{amsfonts}
\usepackage{url}
\usepackage{todonotes}

\newtheorem{lemma}{\bf Lemma}[section]
\newtheorem{theorem}[lemma]{\bf Theorem}
\newtheorem{definition}[lemma]{\bf Definition}
\newtheorem{proposition}[lemma]{\bf Proposition}
\newtheorem{observation}[lemma]{\bf Observation}
\newtheorem{example}[lemma]{\bf Example}
\newcommand{\Pow}{{\mathcal P}}
\newcommand{\proofend}{{\hfill $\Box$}}

\begin{document}

\title{Heritability of K\H{o}nig's Property from finite edge sets}

\author{Marc Kaufmann}
\thanks{The first-named author gratefully acknowledges support by the Swiss National Science Foundation [grant number
200021 192079].}
\address{ETH Z\"urich, CH-8092 Z\"urich,
Switzerland}
\email{marc.kaufmann@inf.ethz.ch}

\author{Dominic van der Zypen}
\address{Bern, Switzerland}
\email{dominic.zypen@gmail.com}


\begin{abstract}
A hypergraph $H = (V,E)$ is said to have {\em K\H{o}nig's Property}
if there is a matching $M\subseteq E$ and $S\subseteq V$ such that 
$|S \cap e| = 1$ for all $e\in M$, and $S$ is a vertex cover of $H$.
Aharoni posed the question whether K\H{o}nig's Property is
inheritable from finite subsets of $E$ (Problem 6.7, \cite{Ah}). We
provide a negative answer and investigate similar questions
for weaker properties.
\end{abstract}
\parindent = 0mm
\parskip = 2 mm
\maketitle


\section{Introduction}

The study of min-max duality in graphs goes back at least to 1912 and the work of Frobenius on determinants of square matrices, a problem which Dénes K\H{o}nig showed to be modelable by the search of perfect matchings in certain bipartite graphs (for an overview of the history cf. Plummer, \cite{Pl}). In 1931, K\H{o}nig proved the first milestone 
\begin{theorem}{(in: \cite{Koe})}
Let $G$ be a finite bipartite graph. Then there exists a vertex cover of its edges whose size is equal to a matching of maximal size.
\end{theorem}
The theorem extended his result for regular bipartite graphs, first communicated seventeen years prior. It was conjectured by Erd\H{o}s that the same result should hold for infinite bipartite graphs as well  and finally settled, in the affirmative, by Aharoni in 1984:
\begin{theorem}{(in: \cite{Ah2})}
In any bipartite graph $G$, there exists a matching $F$ and a cover $C$, such that $C$ consists of the choice of precisely one vertex from each edge in $F$.
\end{theorem}
It is natural to ask when an analogous statement should hold for infinite \textit{hyper}graphs. One approach requires that the result holds for all finite sub\textit{hyper}graphs and transfer the property to the infinite hypergraph of interest. This suffices for graphs (cf. Proposition \ref{finiteineritability}). Aharoni therefore asked if the same heritability holds for hypergraphs. In this article, we give a negative answer to this question.
\section{Preliminaries}

A {\em hypergraph} is a pair $H=(V,E)$
where $V$ is a set and $E\subseteq \Pow(V)$. The elements of
$E$ are called {\em edges}. For all hypergraphs in this
paper we assume that all edges are non-empty (i.e., $\emptyset
\notin E$).

A {\em (vertex) cover} of $H$ is a set
$C\subseteq V$ such that $C \cap e \neq \emptyset$ for all $e\in E$.
The minimal cardinality that a cover can have is said to
be the {\em (vertex) covering number} of $H$, and we denote it by $\nu(H)$.

A {\em matching} is a set $M\subseteq E$ of pairwise disjoint
edges.

In the construction of examples, we will often use $\omega$, the
first infinite ordinal (which is equal to the set of non-negative 
integers).

\begin{observation}\label{obs1}
If $H=(V,E)$ is a hypergraph and $M\subseteq E$ is a matching,
then $|M|\leq \nu(H)$.
\end{observation}

\begin{definition} {\em A hypergraph $H = (V,E)$ is said to have 
{\em K\H{o}nig's Property}
if there is a matching $M_0 \subseteq E$ and a cover $C_0\subseteq V$ such
that $$|C_0 \cap e| = 1 \text{ for all } e\in M_0.$$
}
\end{definition}

The following weakening of K\H{o}nig's Property will be useful:

\begin{definition} {\em A hypergraph $H = (V,E)$  has {\em K\H{o}nig's Weak Property}
if there is a matching $M\subseteq E$ with $|M| = \nu(H)$.
}
\end{definition}

\begin{observation} \label{eqobs} K\H{o}nig's Property implies
K\H{o}nigs Weak Property.
\end{observation}

Note that the contrapositive of \ref{eqobs} says that
if for every matching $M\subseteq E$ we have $|M|<\nu(H)$
then $H$ cannot have K\H{o}nig's Property. We will use
this elementary fact in the construction of the negative
answer to Aharoni's question asked in \cite{Ah}.

The necessary condition of \ref{eqobs}
is not sufficient, as the following example shows:

\begin{example} {\em Let $G$ be $K_\omega$, the complete graph  on $\omega$.
Clearly, $\nu(G) = \aleph_0$, and $$\big\{\{2n, 2n+1\}: n\in \omega\big\}$$
is a matching of cardinality $\aleph_0$. But $G$ does not have K\H{o}nig's 
Property because every vertex cover of $G$ either equals $\omega$ or 
$\omega\setminus\{k\}$ for some $k\in \omega$, and a cover of this form
can never pick exactly $1$ point of any given matching of $G = K_\omega$.
} 
\end{example}

However, for finite graphs the necessary condition of Observation \ref{eqobs} 
is also sufficient:

\begin{proposition} If $G = (V,E)$ is a finite graph with a matching $M\subseteq E$
and $|M| = \nu(G)$, then $G$ satisfies K\H{o}nig's Property.
\end{proposition}
{\em Proof.} Take any maximal matching $M \subseteq E$, a minimal vertex cover $C \subseteq V$ and partition $M$ with respect to $C$:
\[M=(\bigcup_{e \in M, |e \cap C|=1} e)\sqcup (\bigcup_{e \in M, |e \cap C|>1} e) \sqcup (\bigcup_{e \in M, e \cap C= \emptyset} e)\]

(By $\sqcup$ we denote disjoint union.) 

We now prove that only the first union may be non-empty:
\begin{itemize}
    \item $C$ intersects every edge $f \in E$, in particular all $e \in M \subset E$, so \[\bigcup_{e \in M, e \cap C= \emptyset} e = \emptyset\]
    \item Consider $e \in M: |e \cap C|>1$ and any of the $k>1$ vertices $v \in e \cap C$. Either $v$ is contained in another edge of $M$ - contradicting the disjointness of elements of $M$ -  or $C -\{v\}$ yields a vertex cover of $G$. But $|C|$ is finite by assumption, so $\nu(G)=|C|>|C-\{v\}|$, another contradiction. This yields \[\bigcup_{e \in M, |e \cap C|>1} e = \emptyset\]
\end{itemize}
We conclude
\[M=\bigcup_{e \in M, |e \cap C|=1} e\]
\proofend

Another concept related to K\H{o}nig's Property is that of bipartiteness:
\begin{definition} {\em A hypergraph $H = (V,E)$  is {\em bipartite}
if there is $D\subseteq V$ such that 
both $D$ and $V\setminus D$ intersect every $e\in E$ with $|e|>1$.
}
\end{definition}

It is easily seen that any hypergraph satisfying K\H{o}nig's Property
is bipartite.

\section{Heritability from finite edge sets}
In \cite{Ah}, Ron Aharoni asks the following question:
\begin{quote}
{\bf Problem 6.7.} Suppose that every finite
subhypergraph of $H$ (i.e. $H' = (V, E')$
where $E'$ is a finite subset of $E$) satisfies
K\H{o}nig's Property. Does it necessarily follow
that $H$ satisfies K\H{o}nig's Property?
\end{quote}
The following gives a negative answer to this question:
\begin{proposition}\label{prop_one}
Let $E = \{A\subseteq \omega: (\omega \setminus A) \text{ is
finite}\}$. Then
\begin{enumerate}
\item $(\omega, E_0)$ satisfies K\H{o}nig's Property
for all finite $E_0\subseteq E$, but
\item $(\omega, E)$ does not satisfy K\H{o}nig's Property.
\end{enumerate}
\end{proposition}
{\em Proof.}

(1) Let $E_0\subseteq E$ be finite. If $E_0 =\emptyset$, then
the empty cover $C_0 = \emptyset$ and the empty matching $M_0 = \emptyset$
establish K\H{o}nig's Property
in a vacuous way. So we assume that $E_0 \neq \emptyset$.
Since $E_0$ consists of co-finite
subsets of $\omega$, pick $x^* \in \bigcap E_0$. So $C : = \{x^*\}$ is
a cover of $(\omega, E_0)$. Picking any member $e \in E_0$, we see
that the singleton cover $C = \{x^*\}$ together with the
singleton matching $M := \{e\}$ establish K\H{o}nig's Property.

(2) Since every two members of $E$ intersect, we see that
we have $|M| = 1$ for every matching $M\subseteq E$. On the other hand,
$(\omega, E)$ has no finite cover. (If $F\subseteq \omega$ is finite,
then $F$ does not intersect $\omega\setminus F$, which is a member of $E$
by definition.) So $\nu((\omega, E)) = \aleph_0$. So, $(\omega, E)$ 
does not have K\H{o}nig's Weak Property, and therefore 
it does not have  K\H{o}nig's Property (see Observation \ref{eqobs}). \proofend

So we can say that K\H{o}nig's Property is not ``finitely inheritable''.

However, the situation is different if we restrict ourselves
to graphs:

\begin{proposition} \label{finiteineritability}
K\H{o}nig's property is finitely inheritable for graphs.
\end{proposition}
{\em Proof.} Let $G = (V,E)$ be a simple, undirected graph.
If every finite subgraph satisfies K\H{o}nig's Property, then
$G$ has no odd cycles and so $G$ itself is bipartite.
Finally, every bipartite graph satisfies K\H{o}nig's Property,
see \cite{Ah2}. \proofend

In the remainder of this article we investigate the properties
that are weaker than K\H{o}nig's Property for finite inheritability.

We will consider bipartiteness first. 

\begin{example} \label{omegaex} 
{\em Let $[\omega]^\omega$ denote the collection of 
all infinite subsets of $\omega$. We consider the hypergraph $H = 
(\omega, [\omega]^\omega)$.

\begin{enumerate}
    \item for all finite sets $E_0\subseteq [\omega]^\omega$
the hypergraph $(\omega, E_0)$ is bipartite: If $E_0 = \emptyset$, then 
there is a trivial bipartition of 
$(\omega, \emptyset)$. Otherwise, let $$D =
\{ \min(e) : e \in E_0\}.$$ It is easy to 
see that the finite set $D$ intersects every member
of $E_0$, and so does $\omega\setminus D$, as
it is co-finite.
    \item $(\omega, [\omega]^\omega)$ is not bipartite:
    Take any $D \subseteq \omega$. If $D$ is finite,
    then $D$ does not intersect the edge $\omega\setminus D$.
    If $D$ is infinite, then $D$ itself is an edge, which
    gives us the problem that $(\omega\setminus D) \cap D =
    \emptyset$. So, no matter whether $D$ is finite or
    infinite, either $D$ or $\omega\setminus D$ has
    empty intersection with some member of $[\omega]^\omega$.
    Therefore $H = (\omega, [\omega]^\omega)$ is not bipartite.
\end{enumerate}
}
\end{example}
For finite edge size, the situation changes, and with Tychonoff's compactness
theorem \cite{Ty} we can make the following positive statement.

\begin{proposition}\label{compactnessprop}
Let $H= (V,E)$  be a hypergraph with edges of finite size only, such that for 
all finite subsets $E_0\subseteq E$ the hypergraph $(V, E_0)$ is bipartite. 
Then $H$ is bipartite.
\end{proposition}

{\em Proof.} Suppose $H=(V,E)$ is a hypergraph such that all members 
of $E$ are finite, and we assume that $H$ is not bipartite.
We will construct $E_0\subseteq E$ finite such that $(V, E_0)$ is not bipartite.

Let $\{0,1\}^V$ be the set of all maps $f: V \to \{0,1\}$. Then $H$ not 
being bipartite is equivalent to saying that \begin{quote} $(\star)$ 
for every map $f\in \{0,1\}^V$, there is
$e\in E$ with $|e|>1$ and $f$ is constant on $e$. \end{quote}

Let $E_{>1} = \{e\in E: |e|>1\}$. To every $e\in E_{>1}$ we 
associate the set $$C_e = \{f\in \{0,1\}^V: f\text{ is constant on }e\}.$$ 

{\it Claim.} For every $e\in E_{>1}$, the set $C_e$ is open in the product 
topology of $\{0,1\}^V$ (where we endow the base space $\{0,1\}$ with the 
discrete topology). - Fix $e\in E_{>1}$, keeping in mind that $e$ is finite.
For $v\in V$ let $\text{pr}_v: \{0,1\}^V\to \{0,1\}$ be the projection
map sending $f\in \{0,1\}^V$ to $f(v)$.
Then $$C_e = \Big(\bigcap_{v\in V}\text{pr}^{-1}(\{0\})\Big) \cup 
\Big(\bigcap_{v\in V}\text{pr}^{-1}(\{1\})\Big).$$ Since
the intersections involved are finite, $C_e$ is open.

Statement $(\star)$ above amounts to saying that 
$${\mathcal U} := \big\{C_e: e\in E_{>1}\big\}$$
is an (open) cover of $\{0,1\}^V$. Using Tychonoff's theorem, we get a
finite subcover ${\mathcal U}_0\subseteq {\mathcal U}$. 

So there is $n\in\mathbb{N}$ and edges $e_1, \ldots, e_n\in E_{>1}$
such that $${\mathcal U}_0=\{C_{e_1}, C_{e_2},\ldots,C_{e_n}\}.$$

Set $E_0 := \{e_1, e_2,\ldots, e_n\}$.

Recall that every member of $E_0$ has more than $1$ element by definition.
Since ${\mathcal U}_0$ is an open cover of $\{0,1\}^V$, we conclude that 
for every $f\in \{0,1\}^V$ there is $e \in E_0$ such that $f$ is constant
on $e$. This is exactly statement $(\star)$ above, and we conclude
that $(V,E_0)$ is not bipartite. \proofend

Next, we turn to K\H{o}nig's Weak Property.

The example used in \ref{prop_one} shows that for K\H{o}nig's Weak Property is not
finitely inheritable either. But again, things change if we restrict ourselves to 
edges of {\sl finite size}.

\begin{proposition}
Let $H= (V,E)$  be a hypergraph with edges of finite size only, such that for 
all finite subsets $E_0\subseteq E$ the hypergraph $(V, E_0)$ satisfies
K\H{o}nig's Weak Property. Then $H$ satisfies K\H{o}nig's Weak Property.
\end{proposition}

{\em Proof.} First, note that a straightforward application
of Zorn's Lemma implies that there is a matching $M$ that
is maximal with respect to set inclusion $\subseteq$. 
Maximality implies that $\bigcup M$ is a vertex cover. 
So by Observation \ref{obs1} and
definition of $\nu(H)$ we have $$|M| \leq \nu(H) \leq |M|.$$

{\em Case 1.} If $M$ is infinite, then the inequality above
directly implies $|M| = \nu(H)$ because all edges are finite,
so we get K\H{o}nig's Weak Property.

{\em Case 2.} If $M$ is finite, then $\nu(H)$ is finite, and a 
compactness argument similar to the one used in Proposition
\ref{compactnessprop} shows that there is a finite set 
$E_0\subseteq E$ with $\nu((V, E_0)) = \nu(H)$. Since
$(V, E_0)$ is balanced by assumption, there is a matching $M_0
\subseteq E_0$ such that $|M_0| = \nu((V, E_0)) = \nu(H)$, which
finishes the proof. \proofend

\section{Other notions of bipartiteness}
There are several ways to define a bipartite hypergraph. 
(All the definitions agree on graphs.)

\begin{definition} {\em 
Let $H= (V, E)$ be a hypergraph with $|e|>1$ for all $e\in E$. We say 
that $H$ has the {\em choosability
property} (CP) if there is $C \subseteq V$ such that $|C \cap e| = 1$
for all $e\in E$.
}
\end{definition}
Clearly, the choosability property implies bipartiteness. 

For showing that (CP) is not finitely inheritable in general,
we need to introduce the concept of an {\em  
almost disjoint family}.

\begin{definition}{\em
Let $A,B\in [\omega]^\omega$ (the collection of infinite subsets of $\omega$).
We say $A, B$ are {\em almost disjoint} if $A\cap B$ is finite. An 
{\em almost disjoint family} is a set ${\mathcal A} \subseteq [\omega]^\omega$
consisting of pairwise almost disjoint infinite subsets of $\omega$.
}
\end{definition}

A standard application shows that every maximal almost disjoint family
is contained in a maximal almost disjoint (``MAD'') family. 

\begin{example}{\em 
Erd\H{o}s and Shelah constructed a non-bipartite MAD family 
$E\subseteq [\omega]^\omega$, see
\cite{ES}, Theorem 1.1. Let $H = (\omega, E)$.
\begin{enumerate}
    \item Since $H$ is not bipartite, it does not have the choosability property.
    
    \item Let $E_0\subseteq E$ be finite. If $E_0 = \emptyset$, then
    (CP) is satisfied vacuously. So suppose $E_0 \neq \emptyset$. 
    Let $F$ be the set of members of $\omega$ that are
    an element of more than one member of $E_0$. Since $E_0$ is finite, and 
    the intersection of any two distinct members of $E_0$ is finite, 
    we get that $F$ is finite (possibly empty). Let $x_0 = 1 + \max(F)$.
    From every member $e\in E_0$, the set 
    $$e\setminus \underbrace{\{0,\ldots, x_0\}}_{= x_0+1}$$ is non-empty
    because $e$ is infinite. So for every $e\in E_0$, pick one member of
    that set. This establishes the choosability property.
\end{enumerate}
}
\end{example}
Again, by a compactness argument, the choosability property
is finitely inheritable. But the compactness argument is 
quite different from the one used in Proposition \ref{compactnessprop},
so we will state the proof in its entirety.
\begin{proposition}
Let $H= (V,E)$  be a hypergraph with edges of finite size only, such that for 
all finite subsets $E_0\subseteq E$ the hypergraph $(V, E_0)$ has (CP). 
Then $H$ has (CP).
\end{proposition}

{\em Proof.} Consider the Cartesian product $$K=\prod \{e: e\in E\}$$ of all edges, endowed
with the product topology (where every individual edge is given the discrete topology). 
The set $K$ may be understood as the set of simultaneous choices 
$v_e\in e$ for all edges $e$. For every pair $e\ne e'$ we define
$$U(e,e') = \{f\in K: f(e) \in e'\}.$$ It is a standard exercise to
verify that for all $e\neq e' \in E$ the set $U(e,e')$ is open in the
product topology. 

{\em Case 1.} $\bigcup\{U(e,e'): e\neq e'\in E\} \neq K$. It immediately follows
that $H = (V, E)$ has the choosability property.

{\em Case 2.} $\bigcup\{U(e,e'): e\neq e'\in E\} = K$. Then the collection
$\{U(e,e'): e\neq e'\in E\}$ forms an open cover of $K$. So we apply
Tychonoff's compactness theorem and get a finite subcover. In that case,
we obtain a finite $E_0\subseteq E$ such that $(V, E_0)$ does not
have (CP), contradicting the assumption of the Proposition.
\proofend

\end{document}